\documentclass[12pt]{article}

\usepackage[inline]{enumitem}
\usepackage{amsmath}
\usepackage{bbm}
\usepackage{hyperref}
\usepackage{amsfonts}
\usepackage{graphicx}
\usepackage{natbib}
\usepackage{bm}
\usepackage{tikz}
\usepackage[title]{appendix}

\newcommand{\n}[2]{\mathcal{N}\left( #1,#2\right)}
\newcommand{\pr}{\mathcal{P}}
\newcommand{\hpr}{\hat{\mathcal{P}}}
\newcommand{\E}{E}
\newcommand{\Cov}{\mathrm{cov}}
\newcommand{\Var}{\mathrm{var}}

\renewcommand{\baselinestretch}{1}
\usepackage{etoolbox}
\usepackage{authblk}

\title{Demystifying statistical learning based on efficient influence functions}
\date{\today}

\newcommand{\blind}{0}

\begin{document}

\def\spacingset#1{\renewcommand{\baselinestretch}%
{#1}\small\normalsize} \spacingset{1}

%%%%%%%%%%%%%%%%%%%%%%%%%%%%%%%%%%%%%%%%%%%%%%%%%%%%%%%%%%%%%%%%%%%%%%%%%%%%%%

\if0\blind
{
  \title{\bf Demystifying statistical learning based on efficient influence functions}
    \author[1]{Oliver Hines}
    \author[2]{Oliver Dukes}
    \author[1]{Karla Diaz-Ordaz}
    \author[1,2]{Stijn Vansteelandt}
    \affil[1]{Department of Medical Statistics, London School of Hygiene and Tropical Medicine, UK}
    \affil[2]{Department of Applied Mathematics, Computer Science and Statistics, Ghent University, Ghent, Belgium}

  \maketitle
} \fi

\if1\blind
{
  \bigskip
  \bigskip
  \bigskip
  \begin{center}
    {\LARGE\bf Demystifying statistical learning based on efficient influence functions}
\end{center}
  \medskip
} \fi

\bigskip

\begin{abstract}
Evaluation of treatment effects and more general estimands is typically achieved via parametric modelling, which is unsatisfactory since model misspecification is likely. Data-adaptive model building (e.g. statistical/machine learning) is commonly employed to reduce the risk of misspecification. Na\"{\i}ve use of such methods, however, delivers estimators whose bias may shrink too slowly with sample size for inferential methods to perform well, including those based on the bootstrap. Bias arises because standard data-adaptive methods are tuned towards minimal prediction error as opposed to e.g. minimal MSE in the estimator. This may cause excess variability that is difficult to acknowledge, due to the complexity of such strategies.

Building on results from nonparametric statistics, targeted learning and debiased machine learning overcome these problems by constructing estimators using the estimand's efficient influence function under the nonparametric model. These increasingly popular methodologies typically assume that the efficient influence function is given, or that the reader is familiar with its derivation.

In this paper, we focus on derivation of the efficient influence function and explain how it may be used to construct statistical/machine-learning-based estimators. We discuss the requisite conditions for these estimators to perform well and use diverse examples to convey the broad applicability of the theory.

\end{abstract}

\vfill

\spacingset{1.45} % DON'T change the spacing!

\section{Introduction}

The standard statistical approach of building a model, extracting one or more coefficients and reporting their estimates and associated measures of uncertainty (e.g. confidence intervals) is increasingly being criticised (see e.g. \cite{van2015statistics}).
This standard practice encourages the use of overly simplistic, but misspecified models in order to maintain a simple interpretation of the end result \citep{breiman2001statistical}.
It moreover makes the meaning and definition of the reported coefficients dependent upon the selected models. Inference for such `data-dependent' parameters is not straightforward; ignoring their data-dependent nature, as is commonly done, induces bias, excess variability that is not acknowledged by default standard error estimators and, as a result, overly simplistic inferences.

Building on important results on nonparametric estimation of statistical functionals \citep{pfanzagl1985contributions,pfanzagl1990estimation,bickel1993efficient}, \cite{van_der_laan_targeted_2006}, \cite{van_der_laan_targeted_2011}, \cite{Robins2008} and more recently \cite{chernozhukov_double/debiased_2018}, showed how the aforementioned concerns can be accommodated by centering a statistical analysis around a predefined nonparametric estimand. This is a model-free functional of the observed data distribution which characterises the quantity one wishes to infer from data \citep{Berk2021}. It follows from the existing literature that root-$n$ estimators with well understood asymptotic behaviour can often be derived (under feasible conditions) by making use of the estimand's so-called efficient influence function or canonical gradient under the nonparametric model. The resulting strategies are known as `targeted learning' or 'debiased' machine learning, because they effectively enable the use of data-adaptive estimation strategies to model the data-generating distribution, such as variable selection procedures and machine learning algorithms, whilst permitting valid inference of the estimand of interest.

These developments are quite revolutionary in that they are changing the way in which - we believe - data will be analysed in the future. In particular, they shift the focus from model building and validation to choosing estimands that are well connected to scientific questions of interest \citep{Petersen2014}. This shift enables the analysis to be specified before data is obtained, rather than deciding which statistical quantities to report once a model has been validated, as is usually the case e.g. following model/variable selection. Furthermore, model based analyses usually assume the final model was known a priori, whereas estimand inference based on efficient influence functions tend to be `honest' in the sense of expressing also the uncertainty around selecting the data-generating model, see e.g. \cite{Robins2006} for a precise definition of confidence set `honesty'.

The derivation of the efficient influence function is often regarded as somewhat of a `dark art'. One reason is that it is not given much attention in textbooks on the topic and neither is it given much focus in statistics education. Textbooks that refer to such derivations often rely on a fluency in concepts from functional analysis (e.g. Hilbert Spaces). A further reason is that the majority of research articles that derive the efficient influence function of a statistical estimand, rely on manipulating a derivative expression into a canonical form, as the integral of a product of an efficient influence function and a score function. These derivations are often complicated, with some steps appearing as if from nowhere to achieve the desired form.

In this tutorial paper, we instead advocate an equivalent approach based on Gateaux derivatives, formalised by \cite{ichimura2015influence}, which is much simpler in our opinion. We will explain this approach and show how to make use of it, while also providing intuitive insight into what an efficient influence function is. We will moreover explain how root-$n$ converging statistical/machine-learning-based estimators can be constructed, using the efficient influence function, and what conditions are needed for these to work well. This tutorial obeys the principles of van der Laan's `roadmap' \citep{van_der_laan_targeted_2011}. It is aimed to be broadly accessible to students and researchers who would like to derive efficient influence functions for all sorts of nonparametric estimands, using simple differentiation methods, such as the chain rule. We use diverse examples first to show the steps in calculating the efficient influence function (Section \ref{CalcIF}), and also to convey the very broad applicability of the theory (Section \ref{derivations}).

\section{Step 1: Defining the estimand of interest}

The starting point of most statistical analyses is a (semi)parametric model, which is then often interpreted as representing how nature has generated the data.
For certain applications, such as in the physical sciences, this model can be the result of a deep theoretical understanding of the data-generating mechanism. However, oftentimes, especially in the spheres of medicine, psychology and economics, the model is chosen for its simplicity and convenience. Many ubiquitous models, such as the generalized linear and Cox proportional hazards models, are commonly used without reference to a mechanistic understanding, rather because the parameters indexing those models provide useful summaries of associations that are of interest to the analysis. This is problematic for various reasons.
First, nature is rarely as simple as we would like it to be. This leaves many data analysts torn between reporting a simple model, which is likely misspecified, versus reporting a complex model, which is difficult to interpret \citep{breiman2001statistical}.
It demands choosing between an analysis result that is likely biased (as a result of model misspecification) versus one that is likely useless (in view of its complexity).
Second, standard statistical theory for (semi)parametric models was developed for settings where the model is a priori justified by some biological, economic, ... theory (so that one can assume it to be correct)
and where moreover the data analyst commits to using that model. The truth is that a given model is rarely known to be correct, and that data analysts therefore do not commit to a single model, by adopting model selection strategies. This invalidates standard statistical theory.
Third, even the common attempt to infer the model from data (for instance, by relying on variable selection strategies) is overly ambitious as many competing models often fit the data nearly equally well \citep{breiman2001statistical}. While this is generally well realised, it is also then systematically `forgotten' in how we report and interpret statistical analysis results.

To accommodate these concerns, we will instead aim to infer so-called nonparametric estimands. These are functionals of the true observed data distribution $\mathcal{P}$, which are well defined without reference to a (semi)parametric model, and target the scientific question of interest.
With interest in the mean outcome $Y$, such estimand is unambiguously defined as
\[\Psi_1(\mathcal{P})=E_{\mathcal{P}}(Y),\]
where the subscript $\mathcal{P}$ explicates that the expectation $E_{\mathcal{P}}$ is calculated w.r.t. the true distribution $\mathcal{P}$ of $Y$.
We will equivalently write this as
\[\Psi_1(\mathcal{P})=\pr(Y) = \int y d\pr(y). %\label{expectation_1}
\]
where $d\pr(y)$ denotes integration w.r.t. to the probability measure $\pr$ for the random variable $Y$. When $Y$ is continuous, $d\pr(y)$ in this expression can be replaced with $f(y)dy$ to recover the Riemann integral over the probability density function of $Y$. For many of the examples in this paper we work with Riemann integrals, which are likely to be familiar to most readers.

As a second example, suppose we are interested in the effect of a dichotomous exposure $X$ (coded 0 or 1) on an outcome $Y$ in the presence of data on a possibly high-dimensional vector of covariates $Z$ that is sufficient to adjust for confounding. Then a relevant (statistical) estimand could be defined as
\[\Psi_2(\mathcal{P})=E_{\mathcal{P}}\left\{E_{\mathcal{P}}(Y|X=1,Z)-E_{\mathcal{P}}(Y|X=0,Z)\right\},\]
where, with a slight abuse of notation, the subscript $\mathcal{P}$ now explicates that the expectation $E_{\mathcal{P}}$ is calculated w.r.t. the true distribution $\mathcal{P}$ of $(Z,X,Y)$. This is known in the causal inference literature as the \emph{average causal effect} or \emph{average treatment effect}.

Alternatively, regardless of whether the exposure is dichotomous or not, its effect on $Y$ can also be expressed using the estimand
\[\Psi_3(\mathcal{P})=\frac{E_{\mathcal{P}}\left[\left\{X-E_{\mathcal{P}}(X|Z)\right\}Y\right]}{E_{\mathcal{P}}\left[\left\{X-E_{\mathcal{P}}(X|Z)\right\}^2\right]},\]%E_{\mathcal{P}}(Y|A,Z)
which equals the expected conditional covariance between $X$ and $Y$, given $Z$, divided by the expected conditional variance of $X$, given $Z$.
Where we are happy to assume a partially linear model for $E_{\mathcal{P}}(Y|X,Z)$, such as
\[E_{\mathcal{P}}(Y|X,Z)=\beta X+\omega(Z),\]
for some function $\omega(.)$, $\Psi_3(\mathcal{P})$ reduces to $\beta$, but remains well defined outside this model \citep{Robins2008,Vansteelandt2020}.

\section{Step 2: Calculate the estimand's efficient influence function}\label{Sec:IF}

\subsection{Preliminaries}

Throughout, we will assume that we have access to i.i.d.~ observed data on $O_i\equiv (Z_i,X_i,Y_i)$ for subjects $i=1,...,n$.
An estimator of the above estimands is then readily obtained by substituting $\pr$ by an estimator $\hpr_n$, where the sub-index $n$ denotes the sample size.
For instance, choosing $\hpr_n$ to equal the empirical distribution, $\pr_n$, of the observations $Y_1,...,Y_n$ gives rise to the empirical `plug-in' estimator
\[\Psi_1(\hpr_n) = \hpr_n(Y)=\frac{1}{n}\sum_{i=1}^nY_i.\]

For $\Psi_2(\pr)$, let $\hpr_n$ be any distribution of $(Z,X,Y)$ such that the marginal distribution of $Z$ is given by its empirical distribution, and that the conditional distribution of $Y$ given $X=x$ for $x=0,1$ and $Z=Z_i$ for $i=1,...,n$ has conditional mean equal to a given estimator $\hat{E}(Y|X=x,Z=Z_i)$, such as the prediction from some machine learning algorithm. Then
\[\Psi_2(\hpr_n)=\frac{1}{n}\sum_{i=1}^n \hat{E}(Y|X=1,Z=Z_i)-\hat{E}(Y|X=0,Z=Z_i).\]
Finally, for $\Psi_3(\pr)$, let $\hpr_n$ be any distribution of $(Z,X,Y)$ such that the conditional distribution of $X$ given $Z=Z_i$ for $i=1,...,n$ has conditional mean equal to a given estimator $\hat{E}(X|Z=Z_i)$,
and that the marginal distribution of $X-\hat{E}(X|Z)$ and $\left\{X-\hat{E}(X|Z)\right\}Y$ is given by its empirical distribution. Then
\[\Psi_3(\hpr_n)=\frac{\sum_{i=1}^n \left\{X_i-\hat{E}(X|Z=Z_i)\right\}Y_i}{\sum_{i=1}^n \left\{X_i-\hat{E}(X|Z=Z_i)\right\}^2} .\]
The key question now is whether $\Psi(\hpr_n)$ is a good proxy for $\Psi(\pr)$. To understand this, we will scale their difference by $\sqrt{n}$.
When this scaled difference converges in distribution (to a non-degenerate law), then we can roughly say that $\Psi(\hpr_n)$ differs from $\Psi(\pr)$ up to a term of the order 1 over root-$n$. We then say that
$\Psi(\hpr_n)$ converges to $\Psi(\pr)$ at parametric rate, or root-$n$ rate, which is usually the best that we can hope to achieve.

For the sample mean $\Psi_1$, we have that $\hpr_n=\pr_n$ so that this scaled difference equals
\begin{align*}
\sqrt{n}\left\{\Psi_1(\hpr_n) - \Psi_1(\pr)\right\} &= \sqrt{n}(\hpr_n - \pr)Y \\
&= \sqrt{n}\hpr_n(Y-\Psi_1) = \sqrt{n}\pr_n(Y-\Psi_1) \\
&= \frac{1}{\sqrt{n}} \sum_{i=1}^n (Y_i - \mu) \overset{d}{\to} \n{0}{\sigma^2},
\end{align*}
by the classical central limit theorem, where $\mu$ and $\sigma^2$ are the mean and variance of $Y$, respectively.

We are lucky here that the difference $\sqrt{n}\left\{\Psi_1(\hpr_n) - \Psi_1(\pr)\right\}$ can be written in terms of the operator $\sqrt{n}(\pr_n - \pr)$ applied to some $Y$, but this is not generally the case, for the following reasons.
%\footnote{The operator $\mathcal{G}_n = \sqrt{n}(\hpr_n - \pr)$ is known as the empirical process operator. Donsker's theorem tells us what sets of functions, $\mathcal{F}$, we can apply this operator to in order to obtain central limit theorems.}.
First, the difference $\sqrt{n}\left\{\Psi(\hpr_n) - \Psi(\pr)\right\}$ will in general depend on how much $\hpr_n$ differs from $\pr$. This was easy in the above example, where $\hpr_n$ refers to the empirical distribution $\pr_n$ of $Y$, whose behaviour is easy to understand. It is much harder in more general cases where $\hpr_n$ may involve data-adaptive estimators, such as predictions $\hat{E}(Y|X=x,Z=Z_i)$ or $\hat{E}(X|Z=Z_i)$ obtained via machine learning or via parametric model building procedures. For such predictions, we may at best have access to some overall, marginal measure of prediction error, but will often have a poor understanding of the bias and imprecision in these predictions at specific covariate levels $Z_i$.  Second, the difference $\sqrt{n}\left\{\Psi(\hpr_n) - \Psi(\pr)\right\}$ will in general also depend on how sensitive the estimand $\Psi(.)$ is to changes in the data-generating distribution. This is also generally poorly understood given that $\pr$ indexing $\Psi(\pr)$ is an infinite-dimensional parameter (apart from exceptional cases where the observed data is discrete).

The situation thus looks a bit hopeless at this stage, and indeed, we will not succeed to understand the difference $\sqrt{n}\left\{\Psi(\hpr_n) - \Psi(\pr)\right\}$ for arbitrary estimators $\hpr_n$ and arbitrary estimands $\Psi(\pr)$. However, we will see that progress can be made for specific estimators $\hpr_n$, and for estimands $\Psi(\pr)$ that are sufficiently smooth in the data-generating law $\pr$. Before proceeding, we will first formalise the right level of smoothness that is needed.

\subsection{Parametric submodels}\label{sec:param_sub}

To understand how sensitive $\Psi(.)$ is to changes in the data-generating distribution, we will first take a step back.
Rather than examining how $\Psi(.)$ changes as we slightly perturb $\pr$ towards $\hpr_n$, we will study the effect of such perturbation in the direction of a fixed, deterministic distribution, $\widetilde{\pr}$, which, for the purpose of this discussion, we shall assume is absolutely continuous with respect to $\pr$ (i.e. the support of $\widetilde{\pr}$ is contained in the support of $\pr$). %The latter means that observed data values which cannot be observed under $\pr$ (because of the density being zero at those values) also cannot be observed under $\widetilde{\pr}$.
There are many ways in which we may change $\Psi(.)$ to $\widetilde{\pr}$.
Here, we will focus on perturbations in the direction parameterised via the one-dimensional mixture model
\begin{equation}\label{submodel}
\pr_t = t\widetilde{\pr}+(1-t)\pr,
\end{equation}
indexed by $t\in[0,1]$, which is called a \emph{parametric submodel}. This is not a parametric model in the usual sense (given that the true data-generating law $\pr$ is unknown), but is used here as a convenient tool to formalise small perturbations away from $\pr$ in the direction of
$\widetilde{\pr}$. In particular, note that $\pr_0 = \pr$ and $\pr_1 = \widetilde{\pr}$.

The sensitivity of $\Psi(.)$ to changes in the data-generating distribution in the direction of $\widetilde{\pr}$ can now be formalised in terms of the pathwise or directional derivative,
$$\lim _{t\downarrow 0}\left(\frac{\Psi\left(\pr_{t}\right)-\Psi(\pr)}{t}\right) = \frac{d  \Psi(\pr_t)}{d  t}\Big|_{t=0},$$
evaluated at $t=0$,
and in the direction of $\widetilde{\pr}$.
When this limit exists (i.e., is finite), it is called a G\^ateaux derivative. This generalises the concept of a directional derivative to functional analysis, describing how to take the derivative of a function with respect to a function.
Informally, when this derivative `exists' for all regular parametric submodels, then we will say that the estimand is pathwise differentiable.
Here, a `regular' parametric submodel is such that its score $\widetilde{\pr}(O)/\pr(O)-1$ has finite variance, a mild restriction that will be needed to ensure that the derivative $d\pr_t/d  t|_{t=0}$ is (or more precisely, inner products with this score are) well-defined.
In the next paragraph, we will formalise this definition of pathwise differentiability. This formalisation will be practically useful, as it will provide insight what the so-called efficient influence function (also referred to as \textit{canonical gradient}, or \textit{influence curve}) is, how it can be calculated, and why it is useful.

As in \cite{Fisher2020}, we develop some intuition by first considering the special case of discrete data $O$ with support $\{o_1,...,o_k\}$. Then
\[\frac{d \Psi(\pr_t) }{d t}\Big|_{t=0}=\sum_{j=1}^k
\frac{d \Psi(\pr_t)}{d \pr_t(o_j)}\Big|_{t=0}\frac{d \pr_t(o_j)}{d t}\Big|_{t=0}.\]
Here, $d\Psi(\pr_t)/d\pr_t(o_j)$ expresses the estimand's sensitivity to small changes in the observed data law.
The second component expresses how the observed data law changes along the considered path. It is easily verified to equal
\[\frac{d \pr_t(o_j)}{d t}\Big|_{t=0}=\widetilde{\pr}(o_j)-\pr(o_j).\]
The resulting identity
\begin{equation}\label{discrete}
\frac{d\Psi(\pr_t)}{d t}\Big|_{t=0}=\sum_{j=1}^k
\frac{d \Psi(\pr_t)}{d \pr_t(o_j)}\Big|_{t=0}\left\{\widetilde{\pr}(o_j)-\pr(o_j)\right\},\end{equation}
is limiting (by being focussed on discrete data) and ignores that the probabilities $\pr_t(o_1),...,\pr_t(o_k)$ are not variation-independent (i.e., they sum to 1 and thus cannot be changed in arbitrary ways) \citep{Fisher2020}.

We therefore appeal to Riesz's representation theorem,
according to which this derivative, when it exists, can be obtained via integration of a unique `representer' $\phi(O,\pr)$ with finite variance under $\pr$, w.r.t. some measure:

\begin{eqnarray}
\frac{d  \Psi(\pr_t)}{d  t}\Big|_{t=0}
&=& \int \phi(o,\pr)\left\{d\widetilde{\pr}(o)-d\pr(o)\right\}=(\widetilde{\pr}-\pr) \{\phi(O,\pr)\}. \label{IF_Rieszl}
\end{eqnarray}
Contrasting identity \eqref{IF_Rieszl} with \eqref{discrete}, we learn that the representer $\phi(O,\pr)$ is a functional derivative which characterises how sensitive the estimand $\Psi(\pr)$ is to changes in the data-generating distribution $\pr$.
It is referred to as the estimand's canonical gradient, efficient influence curve or efficient influence function (under the nonparametric model).
The existence of a representer with finite variance such that \eqref{IF_Rieszl} holds, essentially expresses that the estimand is sufficiently smooth as a functional of the data-generating law (so that the notion of a `derivative' is well-defined); here, the finite-variance condition expresses that the `derivative' of the estimand w.r.t. the data-generating distribution is finite. Since identity (\ref{IF_Rieszl}) is insensitive to constant, additive shifts in $\phi(O,\pr)$, we will henceforth limit ourselves to mean zero functions (under $\pr$) without loss of generality.

We can now more formally define the estimand to be pathwise differentiable when there exists a mean-zero, finite-variance function $\phi(O,\pr)$ which satisfies \eqref{IF_Rieszl}  for all (regular) parametric submodels. Since the efficient influence function has mean zero, $\pr \{\phi(O,\pr)\}=0$, the derivative of $\Psi(\pr_t)$ w.r.t. $t$ can equivalently be represented as the average of the efficient influence function over the distribution $\widetilde{\pr}$:
\begin{eqnarray}
\frac{d  \Psi(\pr_t)}{d  t}\Big|_{t=0} &=& \widetilde{\pr} \left\{ \phi(O,\pr)\right\} = \E_{\widetilde{\pr}}\left\{\phi(O,\pr)\right\}. \label{Riesz_non_canonical}
\end{eqnarray}
This result forms the basis of how we will calculate
the efficient influence function of an estimand.

\subsection{How to calculate the efficient influence function of an estimand}\label{CalcIF}

There are several ways to derive efficient influence functions.
We here advocate the ``point mass contamination'' strategy that we find simplest. In particular, we will perturb the estimand it in the direction $\widetilde{\pr}$ of a point mass at single observation $\tilde{o}$. Identity \eqref{Riesz_non_canonical} then gives the efficient influence function at observation $o$ directly as
\[\phi(o,\pr)=\frac{d  \Psi(\pr_t)}{d  t}\Big|_{t=0},\]
where the right-hand side is a so-called G\^ateaux derivative. This has the same properties as ordinary derivatives, familiar from calculus, such as the chain rule. This will facilitate calculations. The following examples illustrate this.

Throughout, for convenience, we will implicitly assume that we work with continuous variables, but the results continue to hold for discrete variables, or a mix of discrete and continuous variables, upon swapping sums with integrals, indicators with Dirac delta functions, and probability mass functions with probability density functions, where needed. For our purposes, $\mathbbm{1}_{\tilde{o}}(o)$ denotes the Dirac delta function w.r.t. $\tilde{o}$; i.e., the density of an idealized point mass at $\tilde{o}$, which equals zero everywhere except at $\tilde{o}$ and which integrates to 1.

\paragraph{Example 1 (population mean).} As a first, simple example, consider the  mean of $Y$:
\begin{align*}
\Psi(\pr) &= \pr(Y)=\E_{\pr}(Y)= \int yf(y) dy,
\end{align*}
where $f(y)$ denotes the density function of $Y$ under $\pr$ (which we assume to be absolutely continuous w.r.t. the Lebesgue measure, though results hold more generally).
Perturbing in the direction of a single observation $\tilde{y}$,
\begin{align*}
f_t(y) &
= t\mathbbm{1}_{\tilde{y}}(y) + (1-t)f(y),
\end{align*}
one obtains,
\begin{align*}
\Psi(\pr_t) &= t\int y\mathbbm{1}_{\tilde{y}}(y) dy + (1-t)\E_{\pr}(Y)= t\tilde{y} + (1-t)\Psi(\pr).
\end{align*}
By the chain rule, taking a derivative with respect to $t$ at $t=0$, gives
\begin{align*}
\frac{d  \Psi(\pr_t)}{d  t}\Big|_{t=0} &= \tilde{y} - \Psi(\pr).
\end{align*}
Because this has finite variance, we conclude that $\E(Y)$ is pathwise differentiable with efficient influence function $Y - \Psi(\pr)$. $\Box$

\paragraph{Example 2 (density at a point $y$).} Consider next the density at a given value $y$, $\Psi(\pr) =f(y)$. Under the parametric submodel of Example 1, we readily find that
\begin{align*}
\Psi(\pr_t) &= t\mathbbm{1}_{\tilde{y}}(y) + (1-t)f(y).
\end{align*}
By the chain rule, taking a derivative with respect to $t$ at $t=0$, gives
\begin{align*}
\frac{d  \Psi(\pr_t)}{d  t}\Big|_{t=0} &= \mathbbm{1}_{\tilde{y}}(y)- \Psi(\pr).
\end{align*}
Because the Dirac delta function is unbounded when $Y$ is absolutely continuous w.r.t. Lebesgue measure, and therefore has infinite variance, we conclude that $f(y)$ is not pathwise differentiable.
This lack of smoothness is the result of insufficient information in the data on the density $f(y)$ at the single point $y$. It generally implies that no root-$n$ converging estimators can be constructed.
 $\Box$

\paragraph{Example 3 (average density).} Consider next the average density of $Y$:
\begin{align*}
\Psi(\pr) &= \pr\left\{f(Y)\right\}=\E_{\pr}\left\{f(Y)\right\} = \int f^2(y) dy.
\end{align*}
Under the parametric submodel of Example 1,
\begin{align*}
\Psi(\pr_t) &= \int f_t^2(y) dy.
\end{align*}
By the chain rule, taking a derivative with respect to $t$ at $t=0$, gives
\begin{align*}
\frac{d  \Psi(\pr_t)}{d  t}\Big|_{t=0} &= \int 2f(y)\frac{d}{d  t}f_t(y)\Big|_{t=0}  dy \\
&= 2 \int f(y)\left\{\mathbbm{1}_{\tilde{y}}(y)-f(y)\right\} dy \\
&= 2\left\{f(\tilde{y}) - \Psi(\pr)\right\}.
\end{align*}
Here, we use that the Dirac delta function $\mathbbm{1}_{y}(Y)$ has average
\begin{align*}
\pr \mathbbm{1}_{y}(Y)&= \int \mathbbm{1}_{y}(\tilde{y}) f(\tilde{y})d\tilde{y} = f(y),
\end{align*}
equal to the density at $y$ under the law $\pr$.
Since $2\left\{f(Y) - \Psi(\pr)\right\}$ has finite variance, we conclude that $\E_{\pr}\left\{f(Y)\right\}$ is pathwise differentiable with efficient influence function $2\left\{f(Y) - \Psi(\pr)\right\}$. $\Box$

When perturbing the density $f(o)$ of a vector of observations $O$ in the direction of a point mass at $\tilde{o}$, we have the identity
\[\frac{d f_t(o) }{d  t}\Big|_{t=0} = \mathbbm{1}_{\tilde{o}}(o)-f(o),\]
which was also used in Example 2.
This implies a simple formula for the efficient influence function at $\tilde{o}$ of the estimand $\E_{\pr}\left\{g(O,\pr)\right\}$ for some function $g(O,\pr)$ of $O$ and the true distribution:
\begin{eqnarray}
\frac{d}{d  t}\E_{\pr_t}\left\{g(O,\pr_t)\right\}\Big|_{t=0}
&=&\frac{d}{d  t}\left\{\int g(o,\pr_t) f_t(o) do \right\}\Big|_{t=0}\nonumber\\
&=&\left\{\int \frac{d}{d  t}g(o,\pr_t) f_t(o) do+\int g(o,\pr_t) \frac{d}{d  t}f_t(o) do \right\}\Big|_{t=0}\nonumber\\
&=&E_{\pr}\left\{\frac{d}{d  t}g(o,\pr_t)\right\}\Big|_{t=0}+g(\tilde{o},\pr)- \E_{\pr}\left\{g(O,\pr)\right\}.\label{expectation}
\end{eqnarray}
We apply this general identity in the following example.

\paragraph{Example 4 (covariance).}
The covariance
\[\Psi(\pr)= \E_{\pr}\left[\left\{Y-\E_{\pr}(Y)\right\}\left\{X-\E_{\pr}(X)\right\}\right],\]
can be written as $\E_{\pr}\left\{g(O,\pr)\right\}$ for $O\equiv (X,Y)$ and $g(o,\pr) = \left\{y-\E_{\pr}(Y)\right\}\left\{x-\E_{\pr}(X)\right\}$.
Using \eqref{expectation}, we thus find that
\[\frac{d  \Psi(\pr_t)}{d  t}\Big|_{t=0}  =
 E_{\pr}\left[\frac{d}{d  t}
 \left\{Y-\E_{\pr_t}(Y)\right\}\left\{X-\E_{\pr_t}(X)\right\}\Big|_{t=0}\right]+\left\{\tilde{y}-\E_{\pr}(Y)\right\}\left\{\tilde{x}-\E_{\pr}(X)\right\}-\Psi(\pr).\]
By the chain rule, we further have that
\[\frac{d}{d  t}\left\{Y-\E_{\pr_t}(Y)\right\}\left\{X-\E_{\pr_t}(X)\right\}\Big|_{t=0}
= - \frac{d}{d  t}\E_{\pr_t}(Y)\Big|_{t=0}\left\{X-\E_{\pr}(X)\right\}- \frac{d}{d  t}\E_{\pr_t}(X)\Big|_{t=0}\left\{Y-\E_{\pr}(Y)\right\}.\]
Further applying \eqref{expectation} to $\E_{\pr_t}(Y)$ and $\E_{\pr_t}(X)$, we find that
\[\frac{d}{d  t}\left\{Y-\E_{\pr_t}(Y)\right\}\left\{X-\E_{\pr_t}(X)\right\}\Big|_{t=0}
= - \left\{\tilde{y}-\E_{\pr}(Y)\right\}\left\{X-\E_{\pr}(X)\right\}- \left\{\tilde{x}-\E_{\pr}(X)\right\}\left\{Y-\E_{\pr}(Y)\right\},\]
which has mean zero. Since $\left\{Y-\E_{\pr}(Y)\right\}\left\{X-\E_{\pr}(X)\right\}$ has finite variance, we conclude that the covariance $\Psi(\pr)$ is pathwise differentiable with efficient influence function
\[\left\{Y-\E_{\pr}(Y)\right\}\left\{X-\E_{\pr}(X)\right\}-\Psi(\pr).\] $\Box$

\paragraph{Example 5 (potential outcome mean).}

Let $Y^x$ denote the potential outcome under exposure level $x$, which expresses what value the outcome of a given individual would have taken had his/her exposure been set to $x$ by some intervention. Under the usual identifying assumptions, (positivity, consistency, non interference and conditional exchangeability given $Z$) \citep{hernan2006estimating},
\begin{align*}
\Psi(\pr)&=E_{\pr}\left\{E_{\pr}(Y|X=1,Z)\right\}
\end{align*}
is a statistical estimand of the population mean of $Y^1$.

Perturbing $\pr$ in the direction of a point mass at $(\tilde{z},\tilde{x},\tilde{y})$, we find that
\begin{align*}
\Psi(\pr_t)&= \int y f_t(y|1,z)f_t(z) dydz \\
&= \int y \frac{f_t(y,1,z)f_t(z)}{f_t(1,z)} dydz,
\end{align*}
where $f_t(y|x,z)$ is the conditional density function of $Y$, given $X=x,Z=z$, under the parametric submodel, and  $f_t(y,x,z), f_t(x,z),$ and $f_t(z)$ are the joint density functions of $(Y,X,Z), (X,Z),$ and $Z$, respectively under the parametric submodel. By the chain rule, we thus have that
\begin{align*}
\frac{d  \Psi(\pr_t)}{d  t}\Big|_{t=0} =&  \int y \Bigg\{ \frac{f(z)}{f(1,z)}\frac{d}{d  t} f_t(y,1,z)\Big|_{t=0}- \frac{f(y,1,z)f(z)}{f(1,z)^2} \frac{d}{d  t} f_t(1,z)\Big|_{t=0}+  \frac{f(y,1,z)}{f(1,z)} \frac{d}{d  t} f_t(z)\Big|_{t=0} \Bigg\} dydz\\
 =& \int     y \frac{f(y,1,z)f(z)}{f(1,z)} \Bigg( \frac{\mathbbm{1}_{\tilde{y},\tilde{x},\tilde{z}} (y,1,z)}{f(y,1,z)} - \frac{\mathbbm{1}_{\tilde{x},\tilde{z}} (1,z)}{f(1,z)} + \frac{\mathbbm{1}_{\tilde{z}} (z)}{f(z)} - 1   \Bigg) dydz.
\end{align*}
Evaluating the integral gives the canonical gradient of $\Psi(\pr)$ at $(\tilde{z},\tilde{x},\tilde{y})$:
\begin{align*}
\frac{d  \Psi(\pr_t)}{d  t}\Big|_{t=0} =\frac{\mathbbm{1}_{\tilde{x}}(1)}{\pi(\tilde{z},\pr)} \left\{\tilde{y} - m_1(\tilde{z},\pr)\right\} + m_1(\tilde{z},\pr) - \Psi(\pr),
\end{align*}
where $m_1(z,\pr)=\E_{\pr}(Y|X=1,Z=z)$  and $\pi(z,\pr) = f(1|z) = \E_{\pr}(X|Z=z)$ is the propensity score.
We conclude that $\Psi(\pr)$ is pathwise differentiable with the above efficient influence function.

From this, it readily follows that
$\Psi_2(\pr)$ (the average treatment effect) is pathwise differentiable with the efficient influence function given by
\begin{align*}
\varphi_{1}(O,\pr) - \varphi_{0}(O,\pr) - \Psi_2(\pr)
\end{align*}
where $\varphi_{x}(O,\pr)$ is the `uncentered' efficient influence curve
\begin{align}
\varphi_{x}(O,\pr) = \frac{\mathbbm{1}_{X}(x)}{f(x|Z)} \left\{Y - m(x,Z)\right\} + m(x,Z). \label{AIPW_infl}
\end{align} $\Box$

\paragraph{Example 6 (conditional outcome mean).}
Before moving on to more elaborate examples, we finally consider
\begin{align*}
\Psi(\pr)&=E_{\pr}(Y|X=x),
\end{align*}
for a given value $x$, where $X$ may be (absolutely) continuous (w.r.t. Lebesgue measure).
Perturbing $\pr$ in the direction of a point mass at $(\tilde{x},\tilde{y})$, we find that
\begin{align*}
\Psi(\pr_t)&= \int y \frac{f_t(y,x)}{f_t(x)} dy,
\end{align*}
where $f_t(y,x)$ and $f_t(x)$ are the joint density functions of $(Y,X)$ and $X$, respectively under the parametric submodel. By
the chain rule, we thus have that the canonical gradient is
\begin{align}
\phi(\tilde{o},\pr) = \frac{d  \Psi(\pr_t)}{d  t}\Big|_{t=0} =&  \int y \Bigg\{ \frac{1}{f(x)}\frac{d}{d  t} f_t(y,x)\Big|_{t=0}- \frac{f(y,x)}{f(x)^2} \frac{d}{d  t} f_t(x)\Big|_{t=0}\Bigg\} dy \nonumber \\
 =& \int    \left[ \frac{y}{f(x)} \left\{ \mathbbm{1}_{\tilde{y},\tilde{x}} (y,x)-f(y,x)\right\} - \frac{yf(y,x)}{f(x)^2} \left\{\mathbbm{1}_{\tilde{x}} (x)-f(x)\right\}\right]dy \nonumber \\
 =& \frac{\mathbbm{1}_{\tilde{x}} (x)}{f(x)} \left\{\tilde{y}-E_{\pr}(Y|X=x)\right\}. \label{cond_result}
\end{align}
An issue, however, emerges when one considers the variance of the influence function.
\begin{align*}
\Var\left\{\phi(O,\pr)\right\} &= \int \left(\frac{\mathbbm{1}_{\tilde{x}} (x)}{f(x)}\right)^2 \left\{\tilde{y}-E_{\pr}(Y|X=x)\right\}^2 f(\tilde{y}|\tilde{x}) f(\tilde{x}) d\tilde{y} d\tilde{x} \\
&= \frac{\mathbbm{1}_{x} (x)}{f(x)} \int \left\{\tilde{y}-E_{\pr}(Y|X=x)\right\}^2 f(\tilde{y}|x) d\tilde{y} \\
&= \frac{\mathbbm{1}_{x} (x)}{f(x)} \Var(Y|X=x)
\end{align*}
Since the Dirac delta function $\mathbbm{1}_{x}(x)$ takes an infinitely large value when $X$ is continuous (i.e. when its probability distribution is absolutely continuous w.r.t. Lebesgue measure), we conclude that the conditional mean is not pathwise differentiable in that case.

When $X$ is discrete (as in Example 5), however, then we have that the indicator function $\mathbbm{1}_{x} (x)=1$, so that the variance of the efficient influence function is finite (so long as $\Var(Y|X=x) < \infty$ and $f(x)>0$). $\Box$

The approach that we have adopted in the above examples follows the calculation in \cite{hampel1974influence}.
A second, perhaps more common approach, instead uses the following, canonical form
\begin{eqnarray}
\frac{d  \Psi(\pr_t)}{d  t}\Big|_{t=0} &=& \int \phi(o,\pr)\left\{d\widetilde{\pr}(o)-d\pr(o)\right\}\\
&=& \int \phi(o,\pr)\left(\frac{d\widetilde{\pr}(o)}{d\pr(o)}-1\right)d\pr(o)  \nonumber\\
&=& \int \phi(o,\pr)S(o)d\pr(o)  \nonumber\\
&=&\E_{\pr} \left\{ \phi(O,\pr) S(O)\right\}=\pr \left\{ \phi(O,\pr) S(O)\right\}. \label{IF_Riesz}
\end{eqnarray}
The efficient influence function is then calculated as the unique mean zero function $\phi(O,\pr)$ whose inner product (i.e., covariance) with the score $S(O)$ under a parametric submodel $\pr_t$ equals the pathwise derivative $d  \Psi(\pr_t)/d  t |_{t=0}$, for all parametric submodels; see \cite{Levy2019} for a tutorial. This can be quite laborious, however, since one must manipulate score functions and integral expressions, and moreover solve a functional equation like (\ref{IF_Riesz}) \citep{ichimura2015influence}.

The latter approach nonetheless appears more commonly used because it lends itself easier to semiparametric modelling, where the scores $S(O)$ can now be confined to the scores of those parametric submodels that obey the semiparametric model restrictions.
A further reason for the greater popularity of this approach may be the apparent limitation of the approach advocated in this tutorial, that certain estimands (e.g., example 3) cannot be evaluated at $\pr_t$ because of the use of Dirac delta functions.
\cite{ichimura2015influence} note that this does not invalidate the approach, as it can be resolved by substituting the  Dirac delta function in $\pr_t$ by a probability measure, indexed by a bandwidth $h$, that approaches a point mass when the bandwidth converges to 0. This modification justifies the approach that we adopt, but for simplicity it will be left implicit in the remainder of the work.

\section{Step 3: construct an estimator based on the estimand's efficient influence function}

\subsection{Plug-in bias and how to remove it}

The previous results help us to develop insight into the scaled difference
\begin{equation}\label{scalediff}\sqrt{n}\left\{\Psi(\widetilde{\pr}) - \Psi(\pr)\right\}.\end{equation}
In particular, the canonical gradient gave us a way to express the notion of a functional derivative of the estimand w.r.t. directional changes in the data-generating law. This in turn forms the basis of a functional analog to the Taylor expansion, the so-called von Mises expansion, which is essentially derived from the Taylor series expansion of $\Psi(\pr_t)$ about the point $t=1$ in the one-dimensional parametric submodel.
\begin{eqnarray*}
\Psi(\pr)  &=& \Psi(\widetilde{\pr}) + \frac{d  \Psi(\pr_t)}{d  t}\Big|_{t=1} (0-1) + R(\pr,\widetilde{\pr}),
\end{eqnarray*}
where $R(\pr,\widetilde{\pr})$ is a remainder term of the expansion. This expansion contains the pathwise derivative evaluated at $t=1$, which may be evaluated using an anologue of the Riesz-representation theorem result in \eqref{Riesz_non_canonical},
\begin{eqnarray}
\frac{d  \Psi(\pr_t)}{d  t}\Big|_{t=1} &=& - \pr \left\{ \phi(O,\widetilde{\pr})\right\} = - \E_{\pr}\left\{\phi(O,\widetilde{\pr})\right\}, \label{Riesz_non_canonical_t1}
\end{eqnarray}
details of which are given in the appendix. It follows that the scaled difference of interest, equation \eqref{scalediff}, can be written as
\begin{eqnarray}\label{taylor}
\sqrt{n}\left\{ \Psi(\widetilde{\pr}) - \Psi(\pr) \right\}&=& -\sqrt{n}\pr \left\{ \phi(O,\widetilde{\pr})\right\} -\sqrt{n}R(\pr,\widetilde{\pr})
\end{eqnarray}
where we note that this identity is guaranteed to hold, so long as we impose no restrictions on the remainder term, which we will consider later. Now letting $\widetilde{\pr}$ equal $\hat{\pr}_n$, we thus see that
\begin{align*}
\sqrt{n}\left\{ \Psi(\hat{\pr}_n) - \Psi(\pr) \right\}&= -\sqrt{n}\pr \left\{ \phi(O,\hat{\pr}_n)\right\} -\sqrt{n}R(\pr,\hat{\pr}_n)\\
&\approx -\frac{1}{\sqrt{n}}\sum_{i=1}^n\phi(O_i,\hat{\pr}_n) -\sqrt{n}R(\pr,\hat{\pr}_n),
\end{align*}
Here, the term
\begin{equation}
-\frac{1}{\sqrt{n}}\sum_{i=1}^n\phi(O_i,\hat{\pr}_n) \label{drift}
\end{equation}
does not converge to zero (indeed, it would not even converge to zero if $\pr$ were used in lieu of $\hat{\pr}_n$) and may sometimes even diverge.
This tends not to cause asymptotic bias in $\Psi(\hat{\pr}_n)$ (because the calculation of bias requires further scaling by $1/\sqrt{n}$ and, moreover, $\phi(O,\pr)$ has mean zero and $\pr_n$ is assumed to converge to $\pr$).
However, it biases the scaled difference $\sqrt{n}\left\{ \Psi(\hat{\pr}_n) - \Psi(\pr) \right\}$, thereby invalidating na\"{\i}ve confidence intervals and tests.
Understanding the behaviour of \eqref{drift} is difficult as a result of the non-standard behaviour of statistical/machine learning-based estimators affecting the large sample behaviour of $\hpr_n$, which in turn propagates into the behaviour of $\Psi(\hpr_n)$.
Let us reconsider Example 5 (the potential outcome mean), for instance, where
\begin{align*}
\Psi(\pr)&=E_{\pr}\left\{E_{\pr}(Y|X=1,Z)\right\}.
\end{align*}
A plug-in estimator is readily obtained as
\begin{align*}
\Psi(\hpr_n)&=\frac{1}{n}\sum_{i=1}^n m_1(Z_i,\hpr_n).
\end{align*}
Here, $m_1(z,\hpr_n)$ denotes a data-adaptive estimator of $m_1(z,\pr)=\E_{\pr}(Y|X=1,Z=z)$, e.g. obtained using parametric regression models with variable selection, or via machine learning algorithms.
The plug-in bias term (\ref{drift}) then equals\footnote{Note that we have evaluated the plug-in bias term at the true propensity score because the considered plug-in estimator does not rely on an estimated propensity score. One may alternatively evaluate the plug-in bias term at an estimated propensity score, which will then only affect the remainder term.}
\begin{align*}
&-\frac{1}{\sqrt{n}}\sum_{i=1}^n \frac{\mathbbm{1}_{1}(X_i)}{\pi(Z_i,\pr)} \left\{Y_i - m_1(Z_i,\hpr_n)\right\} + m_1(Z_i,\hpr_n) - \Psi(\hpr_n) =-\frac{1}{\sqrt{n}}\sum_{i=1}^n \frac{\mathbbm{1}_{1}(X_i)}{\pi(Z_i,\pr)} \left\{Y_i - m_1(Z_i,\hpr_n)\right\}\\
&=-\frac{1}{\sqrt{n}}\sum_{i=1}^n \frac{\mathbbm{1}_{1}(X_i)}{\pi(Z_i,\pr)} \left\{Y_i - m_1(Z_i,\pr)\right\}+\frac{1}{\sqrt{n}}\sum_{i=1}^n \frac{\mathbbm{1}_{1}(X_i)}{\pi(Z_i,\pr)} \left\{m_1(Z_i,\pr)-m_1(Z_i,\hpr_n)\right\},
\end{align*}
where the second term will often follow a non-standard distribution. For instance, when $m_1(z,\hpr_n)$ is obtained using parametric regression models with variable selection, it will often follow a mixture distribution for each $z$ as a result of variation in the selected model across repeated samples.

The extent to which the plug-in bias term (\ref{drift}) causes bias is thus generally poorly understood as it inherits the behaviour of $\hat{\pr}_n$, which is complex when data-adaptive methods are used.
Rather than attempting to understand its asymptotic behaviour, a much simpler remedy is therefore to adjust the plug-in estimator in such a way that the this bias is zero. One easy way to do this is by defining a new estimator
\[\Psi(\hat{\pr}_n)+\frac{1}{n}\sum_{i=1}^n\phi(O_i,\hat{\pr}_n)\]
obtained by subtracting an estimate of the plug-in bias from the plug-in estimator. Then, the scaled difference between this so-called one-step estimator and $\Psi(\pr)$ is governed by $-\sqrt{n}R(\pr,\hat{\pr}_n)$, which will generally be much smaller.

We will see later that there are other ways of modifying the plug-in estimator so that the resulting estimator has zero plug-in bias.

\subsection{The von Mises expansion}

In the previous section, we have built some intuition into plug-in bias and how it can be removed. In order to understand the behaviour of the scaled difference between the one-step estimator and $\Psi(\pr)$, a more careful derivation is needed. In particular, because $\pr$ is unknown, we substituted it by the empirical distribution function $\pr_n$ of the observed data, but did not express the error this is adding to the results. Let us therefore take a step back to identity (\ref{taylor}).
By adding and subtracting $\sqrt{n}(\pr_n - \pr) \left\{\phi(O,\pr)\right\}$ and $\sqrt{n}\pr_n \left\{\phi(O,\widetilde{\pr})\right\}$ to the righthand side, we obtain
\begin{eqnarray*}
\sqrt{n}\left\{\Psi(\widetilde{\pr}) - \Psi(\pr)\right\}
&=&\sqrt{n} (\pr_n-\pr)\left\{\phi(O,\pr)\right\}-\sqrt{n}\pr_n\left\{\phi(O,\widetilde{\pr})\right\}\\&&+\sqrt{n} (\pr_n-\pr)\left\{\phi(O,\widetilde{\pr})-\phi(O,\pr)\right\}-\sqrt{n}R(\pr,\widetilde{\pr}).
\end{eqnarray*}
Setting $\widetilde{\pr}$ to $\hpr_n$ one can rewrite the plug-in bias in form of the so-called von Mises expansion:
\begin{eqnarray}
\nonumber \sqrt{n}\left\{\Psi(\hpr_n) - \Psi(\pr)\right\} &=& -\sqrt{n}\pr \left\{ \phi(O,\hpr_n)\right\} -\sqrt{n}R(\pr,\hpr_n)  \\
&=&\frac{1}{\sqrt{n}}\sum_{i=1}^n \phi(O_i,\pr) -\frac{1}{\sqrt{n}}\sum_{i=1}^n\phi(O_i,\hpr_n) \nonumber \\&&+\sqrt{n} (\pr_n-\pr)\left\{\phi(O,\hpr_n)-\phi(O,\pr)\right\}-\sqrt{n}R(\pr,\hpr_n) \label{vonMises}.
\end{eqnarray}
Here, the first term converges to a normal, mean zero variate by the central limit theorem and the unbiasedness of the canonical gradient.
The empirical process term (i.e., the third term in \eqref{vonMises}) and the remainder term $\sqrt{n}R(\pr,\hpr_n)$ can often be shown to converge to zero under conditions that we will come back to.

Since the asymptotic behaviour of $\hpr_n$, and therefore also of the second term, is often poorly understood, popular approaches are designed to remove this drift term from the expansion. This can be done in multiple possible ways.

\paragraph{One-step estimator.}
The first is to rewrite the above expansion as
\begin{eqnarray*}
\sqrt{n}\left\{\Psi(\hpr_n)+\frac{1}{n}\sum_{i=1}^n\phi(O_i,\hpr_n) - \Psi(\pr)\right\}
&=&\frac{1}{\sqrt{n}}\sum_{i=1}^n \phi(O_i,\pr) \\&&+\sqrt{n} (\pr_n-\pr)\left\{\phi(O,\hpr_n)-\phi(O,\pr)\right\}-\sqrt{n}R(\pr,\hpr_n),
\end{eqnarray*}
and thus to calculate the estimator of $\Psi(\pr)$ as the \emph{one-step estimator}
\[\Psi(\hpr_n)+\frac{1}{n}\sum_{i=1}^n\phi(O_i,\hpr_n).\]
In Example 3, this delivers
\[\int f^2(y,\hpr_n)dy+\frac{2}{n}\sum_{i=1}^n\left\{f(y_i,\hpr_n)-\int f^2(y,\hpr_n)dy\right\}=\left\{\frac{2}{n}\sum_{i=1}^n f(y_i,\hpr_n)\right\}-\int f^2(y,\hpr_n)dy.\]
where $f(y,\hpr_n)$ is a density estimator.
For Example 5 we consider two different cases. When the propensity score $\pi(Z_i,\pr)$ is known, for instance in randomized experiments, then one obtains the estimator,
\[\Psi(\hpr_n)+\frac{1}{n}\sum_{i=1}^n \frac{\mathbbm{1}_{1}(X_i)}{\pi(Z_i,\pr)} \left\{Y_i - m_1(Z_i,\hpr_n)\right\} + m_1(Z_i,\hpr_n) - \Psi(\hpr_n),\]
When the propensity score is unknown, as is the case for observational data, it must also be estimated (e.g. using a data-adaptive estimator $\pi(Z_i,\hpr_n)$) and the one-step estimator recovers the augmented IPW estimator,
\[\frac{1}{n}\sum_{i=1}^n \frac{\mathbbm{1}_{1}(X_i)}{\pi(Z_i,\hpr_n)} \left\{Y_i - m_1(Z_i,\hpr_n)\right\} + m_1(Z_i,\hpr_n),\]
This propensity score estimation has consequences for the remainder term; see below.

\paragraph{Estimating equation estimators.}
The second is to force the drift term to be zero by using it as an estimating equation; that is, to calculate an estimator for $\Psi(\pr)$ as the solution to an estimating equation given by this drift term:
\begin{align}
0&=\frac{1}{n}\sum_{i=1}^n\phi(O_i,\hpr_n). \label{eq_drift_term}
\end{align}
This is easy in the above examples, where the efficient influence function is linear in $\Psi(\pr)$. In Example 3, solving the identity
\[0=\frac{2}{n}\sum_{i=1}^n\left\{f(y_i,\hpr_n)-\Psi(\hpr_n)\right\}\]
delivers a different estimator than the one-step estimator, namely
\[\Psi(\hpr_n)=\frac{1}{n}\sum_{i=1}^nf(y_i,\hpr_n),\]
with the advantage that it is guaranteed to be non-negative.
In Example 5, solving the identity
\[0=\frac{1}{n}\sum_{i=1}^n \frac{\mathbbm{1}_{1}(X_i)}{\pi(Z_i,\hpr_n)} \left\{Y_i - m_1(Z_i,\hpr_n)\right\} + m_1(Z_i,\hpr_n) - \Psi(\hpr_n)\]
for $\Psi(\hpr_n)$ delivers the same estimator as the one-step estimator.

\paragraph{Targeted learning.}
The third works instead by tuning the initial estimator $\hpr_n$ such that it forces \eqref{eq_drift_term} to hold, which is the focus of targeted learning approaches \citep{van_der_laan_targeted_2006,van_der_laan_targeted_2011}. For instance, tuning the estimator $\hpr_n$ in Example 5 to a
retargeted estimator $\hpr_n^*$ that satisfies
\[0=\frac{1}{n}\sum_{i=1}^n \frac{\mathbbm{1}_{1}(X_i)}{\pi(Z_i,\hpr^*_n)} \left\{Y_i - m_1(Z_i,\hpr^*_n)\right\},\]
ensures that the one-step estimator reduces to the simple plug-in estimator
\begin{align*}
&\frac{1}{n}\sum_{i=1}^n m_1(Z_i,\hpr^*_n),
\end{align*}
which then has standard asymptotic behaviour. This tuning can be achieved in many ways; for instance, one may leave the propensity score model unchanged by defining $\pi(Z_i,\hpr^*_n)=\pi(Z_i,\hpr_n)$ and tune the outcome model by defining,
\[m_1(Z_i,\hpr^*_n)=m_1(Z_i,\hpr_n)+\hat{\epsilon} \frac{1}{\pi(Z_i,\hpr_n)},\]
where $\hat{\epsilon} $ is chosen to set the plug-in bias to zero, i.e., it is the solution to
\[0=\frac{1}{n}\sum_{i=1}^n \frac{\mathbbm{1}_{1}(X_i)}{\pi(Z_i,\hpr^*_n)} \left\{Y_i - m_1(Z_i,\hpr_n)-\hat{\epsilon} \frac{1}{\pi(Z_i,\hpr_n)}\right\}.\]
Retargeting an initial density estimator in Example 3 is less straightforward because of the difficulty of ensuring that the retargeted density continues to be a proper density.

Under sufficient conditions that ensure the empirical process and remainder terms to converge to zero, it follows from the above expansion that all 3 above approaches deliver an estimator $\Psi(\hpr^*_n)$ whose asymptotic distribution obeys
\begin{equation}\label{adist}
\sqrt{n}\left\{\Psi(\hpr^*_n) - \Psi(\pr)\right\} \overset{d}{\to} \n{0}{\pr \left\{\phi(Y,\pr)^2\right\}}.
\end{equation}
This is a powerful result, which means that the asymptotic efficiency bound for a nonparametric estimand can be derived as the expected square of the efficient influence function. Heuristically, this bound is a nonparametric analogue of the Cramer-Rao lower bound, and estimators of the type in \eqref{adist} are said to be asymptotically efficient, in the sense that they are asymptotically equivalent to the estimator obtained by solving an estimating equation with known rather than estimated influence function:
\[0=\frac{1}{n}\sum_{i=1}^n\phi(O_i,\pr).\]
It thus tells us that the influence function behaves like the score function in parametric estimation \citep{Wasserman2006nonparam}. It motivates why the variance of $\Psi(\hpr_n)$ can be estimated as 1 over $n$ times the sample variance of the efficient influence function (evaluated at $\hpr_n$), without needing to account for the uncertainty in $\hpr_n$.
Identity (\ref{adist}) also motivates why the definition of pathwise differentiability includes the requirement of an efficient influence function with finite variance. Pathwise differentiability of an estimand is therefore tantamount to the existence of (regular) root-$n$ consistent estimators of that estimand.

\subsection{Controlling the empirical process term}

The asymptotic behaviour of the empirical process term
\[\sqrt{n} (\pr_n-\pr)\left\{\phi(O,\hpr_n)-\phi(O,\pr)\right\}\]
is generally difficult to understand when data-adaptive statistical methods are used.
However, it becomes much simpler to understand when the estimator $\hpr_n$ is trained on an independent dataset, as one can then reason conditional on that estimator.
Reasoning as such, a direct application of Chebyshev's inequality shows that the empirical process term converges to zero in probability when the conditional variance of $\phi(O,\hpr_n)-\phi(O,\pr)$, i.e.,
\[\pr\left[\left\{\phi(O,\hpr_n)-\phi(O,\pr)\right\}^2\right]\]
given $\hpr_n$, converges to zero in probability. The latter can often be shown to hold when the estimator $\hpr_n$ converges to $\pr$ in probability (or even weaker conditions that certain functionals of $\hpr_n$ converge to the corresponding functionals of $\pr$ in probability) and certain positivity conditions hold (see for instance \cite{chernozhukov2017double,Vansteelandt2020} for detailed examples). The use of an independent sample in this way is important for shrinking the empirical process term, but contrary to what popular wisdom sometimes seems to suggest, does not eliminate the leading plug-in bias terms on which we have focused.

Because one rarely has independent data available to train $\hpr_n$, \cite{van_der_laan_cross-validated_2011} and \cite{chernozhukov_double/debiased_2018} recommend a cross-fitting procedure, whereby the data is split into $K$ folds.
For each individual $i$ from fold $k=1,...,K$, the efficient influence function for that individual is then evaluated in an estimator $\hpr_n$ trained on the data for all individuals, except those in the $k$-th fold. This usually results in a better asymptotic approximation, as reflected by more accurate standard error estimators obtained as 1 over root-$n$ times the sample standard deviation of those influence functions. However, it may induce some finite-sample bias in the estimator as a result of the data-adaptive estimator $\hpr_n$ being trained on a smaller sample of data.

\subsection{Controlling the remainder term}

To understand the remainder term $\sqrt{n}R(\pr,\hpr_n)$, we return to the von Mises expansion (\ref{vonMises}), from which it is seen to equal
\begin{eqnarray*}
\sqrt{n}R(\pr,\hpr_n)
&=&-\sqrt{n} \pr\left\{\phi(O,\hpr_n)\right\}-\sqrt{n}\left\{\Psi(\hpr_n) - \Psi(\pr)\right\}.
\end{eqnarray*}
In Example 5, this is
\begin{eqnarray*}
\sqrt{n}R(\pr,\hpr_n)
&=&-\sqrt{n}E_{\pr}\left[\frac{\mathbbm{1}_{1}(X)}{\pi(Z,\pr)} \left\{Y - m_1(Z,\hpr_n)\right\} + m_1(Z,\hpr_n) - \Psi(\hpr_n)\right] -\sqrt{n}\left\{\Psi(\hpr_n) - \Psi(\pr)\right\}\\
&=&-\sqrt{n}E_{\pr}\left[\frac{\mathbbm{1}_{1}(X)}{\pi(Z,\pr)} \left\{Y - m_1(Z,\hpr_n)\right\} + m_1(Z,\hpr_n) - \Psi(\pr)\right] \\
&=&-\sqrt{n}E_{\pr}\left[\left\{\frac{\pi(Z,\pr)}{\pi(Z,\pr)} -1\right\}\left\{m_1(Z,\pr) - m_1(Z,\hpr_n)\right\}\right]=0.
\end{eqnarray*}
Hence when the propensity score $\pi(Z,\pr)$ is known, the remainder term is zero. Other estimands are also known to have a zero remainder, such as the average density (see Example 3).

When substituting $\pi(Z,\hpr_n)$ for $\pi(Z,\pr)$,
by the Cauchy-Schwarz inequality, the remainder can be upper bounded by
\[\sqrt{n}E_{\pr}\left[\left\{\frac{\pi(Z,\pr)}{\pi(Z,\hpr_n)} -1\right\}^2\right]^{1/2} E_{\pr}\left[\left\{m_1(Z,\pr) - m_1(Z,\hpr_n)\right\}^2\right]^{1/2}.\]
This converges to zero in probability when $\pi(Z,\hpr_n)$ and $m_1(Z,\hpr_n)$ converge to $\pi(Z,\pr)$ and $m_1(Z,\pr)$, respectively, at faster than $n$ to the quarter rate (and $\pi(Z,\hpr_n)$ is bounded away from zero), which is a typical requirement in the non/semiparametric literature. In this specific example, the remainder also shrinks to zero under more general conditions; $m_1(Z,\hpr_n)$ can be allowed to converge at a slow rate, so long as $\pi(Z,\hpr_n)$ is fast converging. This additional flexibility is sometimes known as `rate double-robustness', and does not apply to remainder terms in general, although it does apply for many common estimands in causal inference/missing data problems \citep{rotnitzky2021characterization}. To obtain fast rates of convergence with flexible methods, we typically rely on strong smoothness/sparsity assumptions (e.g. when $Z$ is high dimensional, $\pi(Z,\pr)$ and/or $m_1(Z,\pr)$ should depend on a small number of the covariates), in addition to well-chosen tuning parameters for the learners.

We refer the reader to \cite{Fisher2020} for a rigorous treatment of the remainder terms of the von Mises expansion, which are usually analysed on a case-by-case basis (see for instance \cite{chernozhukov2017double,Vansteelandt2020} for detailed examples).

\section{Examples}\label{derivations}

In this section, we illustrate the calculation of the canonical gradient for the expected conditional covariance and the average derivative effect, deriving the one-step estimators in both cases. Further examples are provided in Appendix B, which also contains results that readers may find helpful for reference.

\subsection{General results}
\label{standard}

For notational convenience we define an operator, $\partial_t$, applied to an arbitrary function, $g(t)$, as
\begin{align*}
\partial_t g(t) = \frac{d  g(t)}{d  t}\Big|_{t=0}.
\end{align*}
For instance, let $f_t(y,x)$ denote a parametric submodel which disturbs the density $f(y,x)$ of $(Y,X)$ at $(y,x)$ in the direction of a point mass at $(\tilde{y},\tilde{x})$. Then from
\begin{align*}
f_t(y|x) &= \frac{f_t(y,x)}{f_t(x)}
\end{align*}
and using the chain rule and the quotient rule for derivatives, we obtain
\begin{align*}
\partial_t f_t(y|x) &= \partial_t \left\{\frac{f_t(y,x)}{f_t(x)}\right\} \\
&=  \frac{\partial_t f_t(y,x) f(x) - f(y,x)\partial_t f_t(x)}{f^2(x)}\\
&= \frac{1}{f(x)} \left[ \mathbbm{1}_{\tilde{y},\tilde{x}}(y,x)-f(y,x) - \frac{f(y,x)}{f(x)} \left\{\mathbbm{1}_{\tilde{x}}(x)-f(x)\right\}  \right] \\
&= \frac{\mathbbm{1}_{\tilde{x}}(x)}{f(x)} \left\{ \mathbbm{1}_{\tilde{y}}(y)-f(y|x) \right\}.
\end{align*}
Similarly to \eqref{expectation}, this expression may be used to derive the following identity for the conditional expectation of an arbitrary function $g(o,\pr)$, where $o=(y,x')'$:
\begin{align}
\partial_t \E_{\pr_t}\left\{g(O,\pr_t)|X=x\right\} &= \partial_t \int g(o,\pr_t) f_t(y|x) dy \nonumber \\
&= \frac{\mathbbm{1}_{\tilde{x}}(x)}{f(x)} \left[g(\tilde{o},\pr) - \E_{\pr}\left\{g(O,\pr)|X=x\right\} \right] + \E_{\pr}\left\{\partial_t g(O,\pr_t)|X=x\right\}. \label{conditional_expectation}
\end{align}
Such generic expressions are helpful to relate to, as they can be used to speed up derivations. For instance, for the potential outcome mean,
defining $m_1(Z,\pr)=\E_{\pr}(Y|X=1,Z)$,
it readily follows from \eqref{expectation} that
\begin{align*}
\partial_t \E_{\pr_t}\left\{m_1(Z,\pr_t)\right\} &= m_1(\tilde{z},\pr) - \E_{\pr}\left\{m_1(Z,\pr)\right\} + \E_{\pr}\left\{\partial_t m_1(Z,\pr_t)\right\},
\end{align*}
and by \eqref{conditional_expectation}, that
\begin{align*}
\partial_t m_1(Z,\pr_t)\ &=  \frac{\mathbbm{1}_{\tilde{x},\tilde{z}}(1,z)}{f(1,z)} \left\{\tilde{y}-m_1(z,\pr)\right\}+0.
\end{align*}
Averaging over the distribution of $Z$ then delivers
\begin{align*}
\E_{\pr}\left\{\partial_t m_1(Z,\pr_t)\right\} &=  \frac{\mathbbm{1}_{\tilde{x}}(1)}{f(1|z)} \left\{\tilde{y}-m_1(z,\pr)\right\}.
\end{align*}
Hence, we recover the same result as before.

\paragraph{Example 7 (expected conditional covariance).}
Consider the expected conditional covariance,
\begin{align*}
\Psi(\pr) = \E_{\pr}\left[\left\{Y-\E_{\pr}(Y|Z)\right\}\left\{X-\E_{\pr}(X|Z)\right\}\right]
\end{align*}
which appears in hypothesis testing \citep{Shah2018} and in parameter estimation in generalized linear models \citep{Vansteelandt2020}.
Define
\begin{align*}
\Cov_t(Y,X|Z) = \E_{\pr_t}\left[\left\{Y-\E_{\pr_t}(Y|Z)\right\}\left\{X-\E_{\pr_t}(X|Z)\right\}|Z\right].
\end{align*}
Upon noting that $\Psi(\pr) = \E_{\pr}\left\{\Cov(Y,X|Z)\right\} $ is of the form in \eqref{expectation}, we find that
\begin{align*}
\partial_t\Psi(\pr_t) = \Cov(Y,X|\tilde{z}) - \Psi(\pr) + \E_{\pr}\left\{\partial_t \Cov_t(Y,X|Z) \right\}.
\end{align*}
The complication is clearly in the final term, which is of the form in \eqref{conditional_expectation}, hence,
\begin{align*}
\partial_t \Cov_t(Y,X|z) = \frac{\mathbbm{1}_{\tilde{z}}(z)}{f(z)} \left[\left\{\tilde{y}-\E_{\pr}(Y|\tilde{z})\right\}\left\{\tilde{x}-\E_{\pr}(X|\tilde{z})\right\} - \Cov(Y,X|z)\right] \\
+ \E[\partial_t\left\{Y-\E_{\pr_t}(Y|Z)\right\}\left\{X-\E_{\pr_t}(X|Z)\right\}|Z=z].
\end{align*}
Similarly to the covariance example previously, the final term above turns out to be zero. It follows, therefore, that the canonical gradient is
\begin{align*}
\phi(O,\pr) = \partial_t\Psi(\pr_t) = \left\{Y-\E_{\pr}(Y|Z)\right\}\left\{X-\E_{\pr}(X|Z)\right\}- \Psi(\pr),
\end{align*}
and since this has finite variance,  the expected conditional covariance is pathwise differentiable.

Constructing a one-step estimator or estimating equations estimator based on the canonical gradient of the expected conditional covariance is relatively straightforward and in fact both methods will provide the same result in this example. The one-step estimator  takes an original plug-in estimator $\Psi(\hpr_n)$ and adds a correction term
\begin{align*}
\Psi(\hpr_n) + \frac{1}{n} \sum_{i=1}^n \phi(O_i,\hpr_n) &= \frac{1}{n}\sum_{i=1}^n \left\{Y_i-\hat{m}(Z_i)\right\}\left\{X_i-\hat{\pi}(Z_i)\right\}.
\end{align*}
where $\hat{m}(z) = \E_{\hpr_n}(Y|z)$ and $\hat{\pi}(z) = \E_{\hpr_n}(X|z)$. Estimation by this strategy therefore requires additional modelling to obtain the functions $\hat{m}(z)$ and $\hat{\pi}(z)$. $\Box$

\paragraph{Example 8 (average derivative effect).}

This example concerns the average derivative effect estimand of \cite{Hardle1989} with canonical gradient given by \cite{Newey1993}.
We let $m(x,z,\pr) = \E_{\pr}(Y|X=x,Z=z)$ be a conditional response surface which is assumed to be differentiable w.r.t. $x$, with derivative $m^\prime(x,z,\pr)$, and we also introduce a known weight function, $w(x,z)$. The average derivative effect estimand is written
\begin{equation*}
\Psi(\pr) = \E_{\pr}\left\{w(X,Z)m^\prime(X,Z,\pr)\right\}.
\end{equation*}
\cite{Powell1989} showed that, for a differentiable function $g(x,z)$, with derivative w.r.t. $x$, $g^\prime(x,z)$,
\begin{equation*}
\E_{\pr}\left\{w(X,Z)g^\prime(X,Z)\right\} = \E_{\pr}\left\{l(X,Z,\pr)g(X,Z)\right\}
\end{equation*}
under regularity conditions, which require that $X$ is a continuous random variable and that $w(x,z)f(x,z)$ is differentiable w.r.t. $x$ and is zero on the boundary of the support of $X$, where $f(x,z)$ used to denote the joint distribution of $(X,Z)$ under $\pr$. In the above expression, \[l(x,z,\pr) \equiv -w^\prime(x,z) - w(x,z) f^\prime(x,z)/f(x,z),\] and, as before, superscript prime denotes the derivative with respect to $x$. Using \eqref{expectation},
\begin{align*}
\partial_t \Psi(\pr_t) = w(\tilde{x},\tilde{z})m^\prime(\tilde{x},\tilde{z}) - \Psi(\pr) + \E_{\pr}\left\{w(X,Z)\partial_t m^\prime(X,Z,\pr_t)\right\}.
\end{align*}
For the final term, we rely on Powell's identity:
\begin{align*}
\E_{\pr}\left\{w(X,Z)\partial_t m^\prime(X,Z,\pr_t)\right\} &= \partial_t \E_{\pr}\left\{w(X,Z) m^\prime(X,Z,\pr_t)\right\}\\
&= \partial_t \E_{\pr}\left\{l(X,Z,\pr) m(X,Z,\pr_t)\right\}\\
&=  \E_{\pr}\left[l(X,Z,\pr) \frac{\mathbbm{1}_{\tilde{x},\tilde{z}}(X,Z)}{f(X,Z)}\left\{\tilde{y} - m(X,Z,\pr)\right\}\right] \\
&= l(\tilde{x},\tilde{z},\pr)\left\{\tilde{y} - m(\tilde{x},\tilde{z},\pr)\right\}.
\end{align*}
Since this has finite variance, the average derivative effect is pathwise differentiable with canonical gradient
\begin{align*}
\phi(O,\pr) = \partial_t \Psi(\pr_t) = l(X,Z,\pr)\left\{Y - m(X,Z,\pr)\right\} + w(X,Z)m^\prime(X,Z,\pr) - \Psi(\pr).
\end{align*}
Using this efficient influence function, an efficient estimator may be easily derived following the one-step or estimating equation strategy. In this case both will result in the same estimator. Setting the sample average of $\phi(O_i,\hpr_n)$ to zero results in the estimator
\begin{align*}
\Psi(\hpr_n) &= \frac{1}{n}\sum_{i=1}^n l(X_i,Z_i,\hpr_n)\left\{Y_i - m(X_i,Z_i,\hpr_n)\right\} + w(X_i,Z_i)m^\prime(X_i,Z_i,\hpr_n)
\end{align*}
This estimator therefore requires modelling the functions $m(x,z,\pr),m^\prime(x,z,\pr) $ and $l(x,z,\pr)$.

We include some extra examples in the Appendix.
\section{Implementation}

We begin by summarising the steps that need to be followed to go from scientific question to (data-adaptive) estimation described in the previous sections.

\paragraph{Step 1:} Defining the estimand of interest.\\
The estimand $\Psi(\pr)$ is a nonparametrically defined statistical functional which is chosen with reference to the scientific question of interest. The estimand might be motivated for a variety of reasons, such as with reference to causal inference (e.g. example 5), independence testing (e.g. example 7), variable importance (e.g., \cite{Williamson2021}), etc.

\paragraph{Step 2:} Calculating its efficient influence function (under the nonparametric model). There are several ways  to do this.

\begin{enumerate}
    \item  \textit{Point-mass contamination.}
We compute the G\^{a}teaux/ pathwise derivative of $\Psi(\pr)$ at $\pr$ in direction of a probability point mass $\widetilde{\pr}$. We consider a parametric submodel $\pr_t = (1-t) \pr+t \widetilde{\pr}$ for $t\in[0,1]$, which we use to evaluate the efficient influence function,
$$
\phi(o,\pr)=\frac{d\Psi\left( \pr_t \right)}{d t} \Big|_{t=0}.
$$

\item
The most general method is to work from the definition of pathwise differentiability. Define a rich class of submodels $\pr_t$ for $t \in (-\epsilon, \epsilon)$ such that $ \pr_{0} = \pr$ and $S(O)=\frac{d}{d t} \log f_{t}(o)\vert_{t=0}$ is the score function of $t$, where $f_t(o)$ denotes the density/probability point mass of $O$ under $\pr_t$. Next one writes the derivative of $\Psi(\pr_t)$, as the integral,
$$\frac{d\Psi\left( \pr_t \right)}{d t} \Big|_{t=0}=\int \phi(o, \pr) S(o) d \pr(o).$$
By the Riesz representation theorem, $\phi(o,\pr)$ is the efficient influence function.

\item For many new estimands, we can  manipulate expressions so they can be expressed as simpler components (by applying the chain and product rules) and use known influence functions as building blocks.
\end{enumerate}

\paragraph{Step 3:} Obtaining an estimator based on the efficient influence function (such as one-step or TMLE) which admits a first-order representation, provided we use sample splitting or cross-fitting.

We recommend cross-fitting in the following way: First the data is partitioned into $K$ \textit{folds}, i.e. $K$ smaller data sets of (roughly) equal size. Next, for each fold $k$, estimate the nuisance functionals (e.g. via data-adaptive methods) using the rest of all the data, excluding that in fold $k$. Use these nuisance functionals to evaluate the efficient influence function for each observations in fold $k$. After repeating for each fold, one is left with an estimate of the efficient influence function for all observations in the dataset. Finally, use these efficient influence function estimates to evaluate the estimator, and the error in the estimator. See \cite{van_der_laan_cross-validated_2011} and \cite{chernozhukov_double/debiased_2018} for more on cross-fitting.

\subsection{Software}
These steps have already been implemented in a variety of R software packages for several estimands from the causal inference literature, e.g., the   average treatment effect, local average treatment effects (using instrumental variables) and the effect curve for a continuous treatment.

For example, the R packages \texttt{npcausal} and \texttt{DoubleML} \citep{DoubleML2021} both have functions to estimate the ATE, the ATT and LATE (based on instrumental variables), and allow the use of sample splitting and a range of machine learners. These packages also implement other estimands that may be of interest to the biostatistics and econometrics literature.

The targeted learning literature is rich in software package implementations of a variety of causal estimands. A new package implementing previous targeted learning estimators, is  \texttt{tlverse}.
The TMLE estimators implemented as part of \texttt{tlverse} are all fitted using cross-validation by default (for example ‘tmle3’).
TMLE cross-validated estimators of counterfactual means and causal effects are also implemented in the R package \texttt{drtmle}.
The R package \texttt{AIPW} \citep{Zhong2021} implements estimation of the ATE by AIPW (corresponding to the estimating equations or one-step estimator as we have seen) and also a TMLE estimator based on machine learning  algorithms.

Other more advanced causal estimands, such as causal effects for a continuous treatment and optimal treatment rules are implemented in both \texttt{npcausal} or \texttt{tmle3}.

\section{Discussion}

Statistical education still focuses primarily on parametric statistical models, which are assumed to reflect how the data is generated. The inferential theory that is taught does not reflect how data is usually analysed, where models are chosen data-adaptively and different models may fit equally well, especially nowadays, given the increased popularity of machine learning. We therefore believe that many courses would be better focused on translating a scientific question into an nonparametric estimand, and basing inference on its efficient influence function under the nonparametric model.

Courses and textbook treatments on the calculus of influence functions often focus on (semi)parametric models \citep{Tsiatis2006}. The resulting derivations can be challenging, as they require one to respect the restrictions that the model imposes on the observed data distribution. Moreover, they show how one can use these restrictions in order to make efficiency gains. Extracting information from modelling assumptions nevertheless comes at the risk of invalid inference when assumptions are violated. By contrast, our focus is on inference under a nonparametric model. This not only makes the resulting inferences more honest, but can dramatically simplify calculations. Additional efficiency gains are then reserved for special cases when restrictions are known to hold by the study design \citep{zhang2008improving}, or reflect strong pre-existing scientific knowledge \citep{liu2021efficient}.

It is difficult, however, to proceed entirely nonparametrically and avoid regularity conditions all together. Indeed, without assumptions on distribution tails, inference of the mean, Example 1 in the current paper, is impossible \citep{Bahadur1956,Bickel1975}. Likewise, many of our examples rely on working models for statistical functionals, necessitating certain regularity \citep{Robins1997}. For instance, the one-step estimator for Example 5 requires estimating $m_1(Z_i,\hpr_n)$. Whilst flexible data-adaptive/ machine learning estimators can be used, these are better thought of as very highly parametric rather than nonparametric, and make assumptions on the true functional $m_1(z,\pr)$, e.g. that it is smooth in $z$. The crucial difference, however, is that compared with the parametric modelling approach, estimators based on the nonparametric model do not `extract efficiency' from highly parametric modelling assumptions.

Because of the crucial role that efficient influence functions play, we focused on their derivation. Whilst the formal justification of the von Mises expansion relies on concepts from advanced mathematics, calculating the efficient influence function can often be done using techniques covered in a basic calculus course. We have illustrated this for several causal and non-causal statistical functionals (estimands); the method of derivation described can lead to simpler proofs than those in the original research papers.

Influence functions have applications beyond using them to define estimators with zero plug-in bias.
Influence functions capture the stability of estimators to outliers (in fact this is one of their original purposes), which makes them additionally useful to diagnose outliers (as measurements with large influence function values). Recently,  influence functions have started to be used in the machine learning literature too. For example, \cite{Koh2017} used influence functions for interpretability of black-box models, by characterising the impact a data point has on the black-box’s predictions. \cite{curth2020semiparametric} and \cite{kennedy2020optimal} use influence functions as the outcome in machine learning procedures of conditional (e.g. subgroup-specific) estimands.

We hope that our contribution helps demystify the calculation of influence functions and thus encourages their wider adoption.
%T

\bibliography{refs}
\bibliographystyle{apalike}

\section*{Appendix A: Riesz Representation Theorem}

Suppose $\pr$ and $\widetilde{\pr}$ are both absolutely continuous w.r.t. some measure $\nu$ and denote the density functions $f(o) = d\pr(o)/d\nu(o)$ and $\tilde{f}(o) = d\widetilde{\pr}(o)/d\nu(o)$. The density of $\pr_t$ w.r.t. $\nu$ is also well defined
\begin{eqnarray*}
f_t(o) &=& \frac{d\pr_t(o)}{d\nu(o)}
=f(o) + t\left\{\tilde{f}(o) - f(o)\right\}.
\end{eqnarray*}
The score function $S_t(o)$ is the derivative of the log density w.r.t. t
\begin{eqnarray*}
S_t(o) &=& \frac{d \log \{f_t(o)\} }{dt}
=\frac{\tilde{f}(o) - f(o)}{f_t(o)}.
\end{eqnarray*}
It follows that
\begin{eqnarray*}
S_t(o)d\pr_t(o) &=& S_t(o)f_t(o)d\nu(o) \\
&=& d\widetilde{\pr}(o) - d\pr(o).
\end{eqnarray*}
Hence $\pr_t \{S_t(O)\} = 0$. Now we consider the $L_2$ Hilbert space defined using the measure $\pr_t$. This is the set of functions $h(O)$ such that $\pr_t \{h(O)\} = 0, \pr_t \{h(O)^2\} <\infty $ and, letting $g(O)$ be another member of this space we define the inner product $\pr_t\{h(O)g(O)\}$. We refer the interested reader to \cite{Levy2019} for an introduction to these Hilbert spaces. Now, assuming that $d\Psi(\pr_t)/dt$ is a continuous linear functional of $S_t(O)$, which is assumed to be a member of the Hilbert space, we use the Riesz Representation Theorem to obtain
\begin{eqnarray*}
\frac{d  \Psi(\pr_t)}{d  t} &=& \pr_t \left\{ \phi(O,\pr_t) S_t(O)\right\} \\
&=& \int \phi(o,\pr_t)S_t(o)d\pr_t(o)  \\
&=& \int \phi(o,\pr_t)\{d\widetilde{\pr}(o) - d\pr(o)\}  \\
&=&(\widetilde{\pr}-\pr) \left\{ \phi(O,\pr_t)\right\}
\end{eqnarray*}
It follows that this expansion holds for all $t$. Also note that $\pr_t\left\{ \phi(O,\pr_t)\right\} = 0$. In the special cases $t=0$ and $t=1$ this allows us to write
\begin{eqnarray*}
\frac{d  \Psi(\pr_t)}{d t} \Big|_{t=0}  &=& \widetilde{\pr} \left\{ \phi(O,\pr)\right\} \\
\frac{d  \Psi(\pr_t)}{d t} \Big|_{t=1}  &=& - \pr \left\{ \phi(O,\widetilde{\pr})\right\}.
\end{eqnarray*}

\section*{Appendix B}

In this Appendix we derive the following results, which readers might find helpful for reference. Here, $F^{-1}(\tau)$ is the quantile function of $Y$ for known $\tau \in [0,1]$, and $F(y|x)$ is the cumulative distribution function of $Y$ given $X=x$. Also $\Theta(u)$ is a step function which takes the value 1 when $u \geq 0$ and 0 otherwise.
\begin{align*}
\partial_t F_t(y|x) &=  \frac{\mathbbm{1}_{\tilde{x}}(x)}{f(x)} \left\{\Theta(y-\tilde{y})-F(y|\tilde{x})\right\} \\
\partial_t \E_{\pr_t}(Y|Y\leq y) &= \frac{\Theta(y-\tilde{y})}{F(y)} \left\{\tilde{y}- \E_{\pr}(Y|Y\leq y)\right\} \\
\partial_t F_t^{-1}(\tau) &= \frac{\Theta\left\{\tilde{y}-F^{-1}(\tau)\right\} + \tau - 1 }{f\left\{F^{-1}(\tau)\right\}}.
\end{align*}
We also illustrate the steps described in the main paper through two further examples: the interventional direct effect, and the incremental propensity score intervention.

\paragraph{Conditional cumulative distribution function.}
Here we consider the conditional cumulative distribution function, $F(y|x)$ where, $y$ and $x$ are known,
\begin{align*}
F(y|x) &= E_{\mathcal{P}}\left\{ \Theta(Y-y)|X=x \right\}
\end{align*}
It is fairly straightforward to recycle the result in \eqref{cond_result}, with $Y$ replaced with $\Theta(Y-y)$ to recover the desired form.

\paragraph{Tail conditional expectation.}
Here we consider the tail conditional expectation, $\E_{\pr}(Y|Y\leq y)$, where $y$ is known:
\begin{align*}
\E_{\pr}(Y|Y\leq y) =\frac{\E_{\pr}\{\Theta(y-Y)Y\} }{F(y)}.
\end{align*}
Now perturbing in the direction of the parametric submodel, and applying the quotient rule,
\begin{align*}
\partial_t \E_{\pr_t}(Y|Y\leq y) &= \frac{\partial_t\E_{\pr_t}\{\Theta(y-Y)Y\}F(y) - \E_{\pr}\{\Theta(y-Y)Y\}\partial_tF_t(y) }{F(y)^2} \\
&= \frac{[\Theta(y-\tilde{y})\tilde{y} - \E_{\pr}\{\Theta(y-Y)Y\}]F(y) - \E_{\pr}\{\Theta(y-Y)Y\}\{\Theta(y-\tilde{y}) - F(y)\} }{F(y)^2} \\
&= \frac{\Theta(y-\tilde{y})}{F(y)} \{\tilde{y}- \E_{\pr}(Y|Y\leq y)\}.
\end{align*}
We notice that the resultant efficient influence function is zero for observations where $\tilde{y}>y$. This coheres with our intuition that the distribution of $Y$ outside the region $Y\leq y$ does not contribute to the asymptotic efficiency bound of $\E(Y|Y\leq y)$.

\paragraph{Quantile function.}
Here we consider the quantile function, $F_t^{-1}(\tau)$, of a continuous random variable $Y$, where $\tau \in [0,1]$ is known. An alternative derivation of the influence curve can be found in \cite{Vaart2013}. We define the estimand $\Psi = \Psi_\tau(\pr) = F^{-1}(\tau)$. The distribution quantile is implicitly defined by
\begin{align*}
\int_{a}^{\Psi_\tau(\pr)} f(y)dy = \tau,
\end{align*}
where $a$ denotes the lower boundary of the support of $Y$ and $f(y)$ is the density function of $Y$. Under the parametric submodel,
\begin{align*}
\int_{a}^{\Psi_\tau(\pr_t)} f_t(y)dy = \tau.
\end{align*}
Differentiating both sides with respect to $t$, the Leibniz integral rule gives us that
\begin{align*}
f_t\left\{\Psi_\tau(\mathcal{P}_t)\right\} \frac{d \Psi_\tau(\mathcal{P}_t)}{d t} + \int_{a}^{\Psi_\tau(\mathcal{P}_t)} \frac{d f_t(y)}{d t}dy = 0.
\end{align*}
Hence,
\begin{align*}
\partial_t \Psi_\tau(\pr_t)  &= \frac{-1}{f\left\{\Psi(\pr)\right\}}\int_{a}^{\Psi(\pr)} \left\{\mathbbm{1}_{\tilde{y}}(y)-f(y)\right\}dy \\
&= \frac{1}{f\left\{\Psi(\pr)\right\}} \left\{ \int_{a}^{\Psi(\pr)} f(y)dy -  \int_{a}^{\Psi(\pr)} \mathbbm{1}_{\tilde{y}}dy \right\} \\
&= \frac{\tau - \left[1 - \Theta\left\{\tilde{y}-\Psi(\pr)\right\}\right] }{f\left\{\Psi(\pr)\right\}}.
\end{align*}
The resulting efficient influence function can be rewritten by defining the function $\rho^\prime_\tau(u) = \Theta(u) + \tau - 1$, which is the derivative (almost everywhere) of the standard quantile regression loss function, $\rho_\tau(u) = u[\Theta(u) + \tau - 1]$. Doing so results in
\begin{align*}
\phi(y,\pr) &= \rho^\prime_\tau\left\{y-\Psi(\pr)\right\}/{f\left\{\Psi(\pr)\right\}}.
\end{align*}
Interestingly, and as an aside, one might wonder how this estimand behaves for different distributions. Let's consider the median when $Y$ follows a univariate normal distribution with mean $\mu$ and standard deviation $\sigma$. For the normal distribution the mean is equal to the median, so $\Psi(\pr) = \mu$. And hence
\begin{align*}
\phi(y,\pr) &= \sigma\sqrt{2\pi}\rho^\prime_{1/2}(y-\mu) \\
\phi(y,\pr)^2 &= \frac{\pi}{2}\sigma^2.
\end{align*}
The standard error in the median estimator is therefore $\E\left\{\phi(Y,\pr)^2/n\right\}^{1/2} \approx 1.253\frac{\sigma}{\sqrt{n}}$. This is 25\% larger than the standard error in the sample mean, which (under the assumption of normality) estimates the same quantity, but achieves the Cramer-Rao lower bound.

\paragraph{Example 9 (interventional direct effect).}

In this example we will derive (one half of) the efficient influence function for the interventional direct effect for mediation, first defined by \cite{Vansteelandt2017}, with an efficient influence function given in \cite{Benkeser2020}.
This estimand is derived using a causal framework and is used to evaluate the effect of a binary outcome, $X$, on an outcome, $Y$, through a set of mediating variables, $M$, given a set of confounder variables, $Z$. Under standard causal assumptions the estimand may be written as a functional of the observed data. We shall not detail these assumptions here, since, once a functional of the data generating distribution is obtained, the causal assumptions are no longer required to derive estimators and efficiency results for it. For our purposes, it is sufficient to define the estimand over the set of variables, $O=(Y,M,X,Z)$, with conditional response surface, $b(m,x,z) = \E(Y|M=m,X=x,Z=z)$,
\begin{equation}
\Psi(\pr) = \int b(m,x^1,z)f(m|x^0,z)f(z) dm dz, \label{interventional_effect}
\end{equation}
where $x^1$ and $x^0$ are known values. Under the parametric submodel,
\begin{equation}
\Psi(\pr_t) = \int b_t(m,x^1,z)f_t(m|x^0,z)f_t(z) dm dz.
\end{equation}
Applying the derivative operator gives
\begin{align*}
\partial_t \Psi(\pr_t) &= \int  \Bigg[ & \partial_t b_t(m,x^1,z) f(m|x^0,z)f(z)\\ &+ b(m,x^1,z) \partial_t f_t(m|x^0,z) f(z) \\ &+  b(m,x^1,z)  f(m|x^0,z) \partial_t f_t(z) \Bigg] dm dz.
\end{align*}
Evaluating these derivatives gives
\begin{align*}
\partial_t \Psi(\pr_t) = \int \Bigg[ & \frac{\mathbbm{1}_{\tilde{o}}(m,x^1,z)}{f(m,x^1,z)}\left\{\tilde{y} - b(m,x^1,z)\right\} f(m|x^0,z)f(z) \\
                                            &+ b(m,x^1,z)\frac{\mathbbm{1}_{\tilde{o}}(x^0,z)}{f(x^0,z)}\left\{\mathbbm{1}_{\tilde{m}}(m)-f(m|x^0,z)\right\} f(z)  \\
                                            & +b(m,x^1,z)  f(m|x^0,z) \left\{\mathbbm{1}_{\tilde{z}}(z)-f(z)\right\}\Bigg]  dm dz
\end{align*}
and evaluating the integral results in the efficient influence function
\begin{align*}
\frac{\mathbbm{1}_{x^1}(X)f(M|x^0,Z)}{f(M,x^1|Z)} \left\{Y - b(M,x^1,Z)\right\} +  \frac{\mathbbm{1}_{x^0}(X)}{f(x^0|Z)} \{b(M,x^1,Z)-a(x^1,x^0,Z)\} + a(x^1,x^0,Z) - \Psi(\mathcal{P}),
\end{align*}
where we define
\begin{align*}
a(x^1,x^0,z) = \int b(m,x^1,z)f(m|x^0,z)dm.
\end{align*}

\paragraph{Example 10 (incremental propensity score intervention).}

The incremental propensity score intervention estimand is motivated by, and derived in the work of \cite{Kennedy2019}. It is an interesting example, since it uses a stochastic intervention which is a function of the true data generating distribution. We define the estimand over the set of variables $O=(Y,X,Z)$, where $X$ is binary with propensity score $\pi(z)=\E_{\pr}(X|Z=z)$, and conditional response surface, $m(x,z)= \E(Y|X=x,Z=z)$,
\begin{equation*}
\Psi(\pr) = \sum_{x=0}^1 \int m(x,z)g_{\pr}(x|z)f(z)dz,
\end{equation*}
where $g_{\pr}(x|z)$ is a probability mass function, which is dependent on the true data generating distribution. \cite{Kennedy2019} propose the `propensity score intervention' indexed by a known value $\epsilon$,
\begin{equation*}
g_{\pr}(x|z) = \frac{x\epsilon \pi(z) + (1-x)\left\{1-\pi(z)\right\}}{\epsilon \pi(z) + 1 - \pi(z)}.
\end{equation*}
This propensity score intervention is motivated by a multiplication on the odds ratio scale,
\begin{equation*}
\frac{g_{\pr}(1|z)}{g_{\pr}(0|z)} = \epsilon \frac{\pi(z)}{1-\pi(z)}
\end{equation*}
although for the purposes of influence function derivation, we are not too concerned with interpretation of the estimand. Under the parametric submodel,
\begin{equation*}
\Psi(\pr_t) = \sum_{x=0}^1 \int m_t(x,z)g_{\pr_t}(x|z)f_t(z)dz.
\end{equation*}
Applying the $\partial_t$ operator gives
\begin{align*}
\partial_t \Psi(\pr_t) = \sum_{x=0}^1 \int  \Bigg[ & \frac{\mathbbm{1}_{\tilde{o}}(x,z)}{f(x,z)} \left\{\tilde{y} - m(x,z)\right\} g_{\pr}(x|z)f(z)
+  m(x,z)\frac{d g_{\pr}(x|z)}{d \pi } \frac{\mathbbm{1}_{\tilde{z}}(z)}{f(z)} \left\{\tilde{x} - \pi(z)\right\}f(z) \\
+ & m(x,z)g_{\pr}(x|z) \left\{\mathbbm{1}_{\tilde{z}}(z) - f(z)\right\} \Bigg] dz,
\end{align*}

where
\begin{align*}
\frac{d g_{\pr}(x|z)}{d \pi } &= \frac{(2x-1)\epsilon}{(\epsilon \pi(z) + 1 - \pi(z))^2} \\
&=\frac{g_{\pr}(1|z) g_{\pr}(0|z)}{ \pi(z)(1- \pi(z))} \left\{\mathbbm{1}_{1}(x)-\mathbbm{1}_{0}(x)\right\}.
\end{align*}

Now, integrating over $z$, becomes
\begin{align*}
\partial_t \Psi(\pr_t) = \sum_{x=0}^1 & \Bigg[\frac{\mathbbm{1}_{\tilde{x}}(x)}{f(x|\tilde{z})} \left\{\tilde{y} - m(\tilde{x},\tilde{z})\right\} g_{\pr}(x|\tilde{z}) \\
+ &m(x,\tilde{z}) \frac{g_{\pr}(1|\tilde{z}) g_{\pr}(0|\tilde{z})}{ \pi(\tilde{z})\left\{1- \pi(\tilde{z})\right\}} \left\{\mathbbm{1}_{1}(x)-\mathbbm{1}_{0}(x)\right\}\left\{\tilde{x} - \pi(\tilde{z})\right\}
+ m(x,\tilde{z})g_{\pr}(x|\tilde{z}) \Bigg] - \Psi(\pr).
\end{align*}
Performing the summation over $x$, the efficient influence function becomes
\begin{align*}
g_{\pr}(1|Z)\varphi_{1}(O,\pr) + g_{\pr}(0|Z)\varphi_{0}(O,\pr)  + \frac{g_{\pr}(1|Z) g_{\pr}(0|Z)}{ \pi(Z)\left\{1- \pi(Z)\right\}}\left\{X - \pi(Z)\right\}\left\{m(1,Z)-m(0,Z)\right\} -\Psi(\pr),
\end{align*}
where $\varphi_{x}(O,\pr)$ is the `uncentered' AIPW influence function as in \eqref{AIPW_infl}.

\end{document}